# Inequalities for Integer and Fractional Parts

**Mihály Bencze**
Department of Mathematics, Áprily Lajos College, Brasov, Romania

**Florentin Smarandache**
Chair of Math & Science Dept., University of New Mexico, Gallup, NM 87301, USA

**Abstract**: In this paper we present some new inequalities relative to integer and functional parts.

**Theorem 1.** If $x > 0$, then $\dfrac{[x]}{3x + \{x\}} + \dfrac{\{x\}}{3x + [x]} \geq \dfrac{4}{15}$, where $[\cdot]$ and $\{\cdot\}$ denote the integer part, and respectively the fractional part.

**Proof.** In inequality $\dfrac{a}{a + 2b + 2c} + \dfrac{b}{2a + b + 2c} + \dfrac{c}{2a + 2b + c} \geq \dfrac{3}{5}$, we take $a = x$, $b = [x]$, $c = \{x\}$.

**Theorem 2.** If $a, b, c, x > 0$, then
$$\dfrac{a}{[x]b + \{x\}c} + \dfrac{b}{[x]c + \{x\}a} + \dfrac{c}{[x]a + \{x\}b} \geq \dfrac{3}{x}.$$

**Proof.** In inequality $\dfrac{a}{ub + vc} + \dfrac{b}{uc + va} + \dfrac{c}{ua + vb} \geq \dfrac{3}{u + v}$, we take $u = [x]$ and $v = \{x\}$.

**Theorem 3.** If $x > 0$ and $a \geq 1$, then
$$\dfrac{[x]}{(a+1)[x] + 2\{x\}} + \dfrac{[x]}{(a+1)\{x\} + 2[x]} \leq \dfrac{2a + 1}{(a+1)(a+2)}.$$

**Proof.** In inequality $\dfrac{x}{ax + y + z} + \dfrac{y}{x + ay + z} + \dfrac{z}{x + y + az} \leq \dfrac{3}{a + 2}$, we take $y = [x]$ and $z = \{x\}$.

**Theorem 4.** If $x > 0$, then

$$[x]\left(\frac{1}{x[x]+x+1}+\frac{1}{[x]\{x\}+[x]+1}\right)+\{x\}\left(\frac{1}{x[x]+x+1}+\frac{1}{x\{x\}+\{x\}+1}\right)\leq 1.$$

**Proof.** In inequality $\dfrac{x}{xy+x+1}+\dfrac{y}{yz+y+1}+\dfrac{z}{zx+z+1}\leq 1$, we take $y=[x]$ and $z=\{x\}$.

**Theorem 5.** If $x>0$, then

$$\frac{x^3}{[x]\left(3[x]^2+3[x]\{x\}+\{x\}^2\right)}+\frac{x[x]^2}{\{x\}\left([x]^2+[x]\{x\}+\{x\}^2\right)}+\frac{x\{x\}^2}{[x]^2+3[x]\{x\}+3\{x\}^2}\geq\frac{3}{2}$$

**Proof.** In inequality $\sum\dfrac{x^2}{y(x^2+xy+y^2)}\geq\dfrac{3}{x+y+z}$, we take $y=[x]$ and $z=\{x\}$.

**Theorem 6.** If $x>0$, $\dfrac{1}{[x]+2\{x\}}+\dfrac{1}{2[x]+\{x\}}\geq\dfrac{1}{x}$.

**Proof.** In inequality $\sum\dfrac{a^2+bc}{b+c}\geq a+b+c$, we take $a=x$, $b=[x]$, $c=\{x\}$.

**Theorem 7.** If $x>0$,

$$\frac{[x]^3}{[x]^2+[x]\{x\}+\{x\}^2}+\frac{\{x\}^3}{3\{x\}^2+3[x]\{x\}+[x]^2}\geq\frac{x\left(3[x]^2-\{x\}^2\right)}{3\left(3[x]^2+3[x]\{x\}+\{x\}^2\right)}$$

**Proof.** In inequality $\sum\dfrac{a^3}{a^2+ab+b^2}\geq\dfrac{a+b+c}{3}$, we take $a=x$, $b=[x]$, $c=\{x\}$.

**Theorem 8.** If $x>0$, then

$$\frac{1}{2[x]^3+4[x]^2\{x\}+4[x]\{x\}^2+\{x\}^3}+\frac{1}{[x]^3+[x]^2\{x\}+[x]\{x\}^2+\{x\}^3}+$$

$$+\frac{1}{[x]^3+4[x]^2\{x\}+4[x]\{x\}^2+2\{x\}^3}\leq\frac{1}{x[x]\{x\}}.$$

**Proof.** In inequality $\sum\dfrac{1}{a^3+b^3+abc}\leq\dfrac{1}{abc}$, we take $a=x$, $b=[x]$, $c=\{x\}$.

**Theorem 9.** If $x>1$, then $4\left(\dfrac{[x]^3}{\{x\}}+\dfrac{\{x\}^3}{[x]}\right)\geq[x]^2+[x]\{x\}+\{x\}^2$.

**Proof.** In inequality $\sum\dfrac{1}{a}(-a+b+c)^3\geq a^2+b^2+c^2$, we take $a=x$, $b=[x]$, $c=\{x\}$.

**Theorem 10.** If $x > 0$, then
$$\frac{x^4}{[x]^2 - [x]\{x\} + \{x\}^2} + \frac{x([x]^3 + \{x\}^3)}{[x]^2 + [x]\{x\} + \{x\}^2} \geq \frac{3}{2}(x^2 + [x]\{x\})$$

**Proof.** In inequality $\sum \frac{a^3}{b^2 - bc + c^2} \geq \frac{3\sum ab}{\sum a}$, we take $a = x$, $b = [x]$, $c = \{x\}$.

**Theorem 11.** If $x > 0$, then $\left|\frac{[x]\{x\}([x] - \{x\})}{x(x + [x])(x + \{x\})}\right| < 1$.

**Proof.** In inequality $\left|\sum \frac{a-b}{a+b}\right| < 1$ we take $a = x$, $b = [x]$, $c = \{x\}$.

**Theorem 12.** If $x > 0$, then $\sqrt{\frac{[x]}{x + \{x\}}} + \sqrt{\frac{\{x\}}{x + [x]}} > 1$.

**Proof.** In inequality $\sum \sqrt{\frac{x}{y+z}} > 2$, we take $y = [x]$, $z = \{x\}$.

**Theorem 13.** If $x > 1$, then $3 + \frac{\{x\}}{[x]} + \frac{[x]}{\{x\}} \geq 3\sqrt[3]{\frac{(x+[x])(x+\{x\})}{[x]\{x\}}}$.

**Proof.** In inequality $\left(\sum a\right)\left(\sum \frac{1}{a}\right) \geq 3\left(1 + \sqrt[3]{\frac{\prod(a+b)}{abc}}\right)$, we take $a = x$, $b = [x]$, $c = \{x\}$.

**Theorem 14.** If $x > 0$, then $\left(\sqrt{[x]} + \sqrt{\frac{[x]\{x\}}{x}} + \sqrt{\{x\}}\right)^4 \geq 32[x]\{x\}$.

**Proof.** In inequality $\sum \sqrt{xy} \geq 2\sqrt[4]{xyz\sum x}$, we take $y = [x]$, $z = \{x\}$.

**Theorem 15.** If $x > 0$, then $(x^2 + [x]\{x\})^2 \geq 6x^2[x]\{x\}$.

**Proof.** In inequality $\left(\sum xy\right)^2 \geq 3xyz\sum x$, we take $y = [x]$, $z = \{x\}$.

**Theorem 16.** If $x > 0$, then
$$x^2 - x\sqrt{[x]\{x\}} + [x]\{x\} \geq \left([x]\sqrt{\{x\}} + \{x\}\sqrt{[x]}\right)\sqrt{x}.$$

**Proof.** In inequality $\sum xy \geq \sum x \sqrt{yz}$, we take $y = [x]$, $z = \{x\}$.

**Theorem 17.** If $x > 0$, then $\sqrt{[x](x + \{x\})} + \sqrt{\{x\}(x + [x])} \leq (2\sqrt{2} - 1)x$.

**Proof.** In inequality $\sum \sqrt{x(y+z)} \le \sqrt{2}\sum x$, we take $y = [x], z = \{x\}$.

**Theorem 18.** If $x > 0$, then $\dfrac{[x]}{x+\{x\}} + \dfrac{\{x\}}{x+[x]} \ge \dfrac{1}{2}$.

**Proof.** In inequality $\sum \dfrac{a}{b+c} \ge \dfrac{3}{2}$, we take $a = x, \ b = [x], \ c = \{x\}$.

**Theorem 19.** If $x > 0$, then $(x+[x])^3 + (x+\{x\})^3 \ge 21x[x]\{x\} + [x]^3 + \{x\}^3$.

**Proof.** In inequality $\sum(x+y)^3 \ge 21xyz + \sum x^3$, we take $y = [x], z = \{x\}$.

**Theorem 20.** If $x > 1$, then $\sqrt{\dfrac{x+[x]}{x+\{x\}}} + \sqrt{\dfrac{x+\{x\}}{x+[x]}} \le \sqrt{\dfrac{[x]}{\{x\}}} + \sqrt{\dfrac{\{x\}}{[x]}}$.

**Proof.** In inequality $\sqrt{\dfrac{x+y}{x+z}} + \sqrt{\dfrac{x+z}{x+y}} \le \dfrac{y+z}{\sqrt{yz}}$, we take $y = [x], z = \{x\}$.

**Theorem 21.** If $x > 0$, then $\dfrac{x}{x+[x]} + \dfrac{x}{x+\{x\}} \ge \dfrac{5}{2}$.

**Proof.** In inequality $2\sum \dfrac{1}{x+y} \ge \dfrac{9}{\sum x}$, we take $y = [x], z = \{x\}$.

**Theorem 22.** If $x > 1$, then $\dfrac{\{x\}}{[x]} + \left(\dfrac{\{x\}}{[x]}\right)^2 + \left(\dfrac{[x]}{\{x\}}\right)^2 + \dfrac{\{x\}^2}{x^2} \ge \left(\dfrac{1}{x} + \dfrac{1}{\{x\}}\right)[x]$.

**Proof.** In inequality $\sum \dfrac{x^2}{y^2} \ge \sum \dfrac{x}{z}$, we take $y = [x], z = \{x\}$.

**Theorem 23.** If $x > 0$, then
$$[x]^2 - [x]\{x\} + \{x\}^2 \ge \dfrac{3}{4}\max\left\{[x]^2; ([x]-\{x\})^2; \{x\}^2\right\}.$$

**Proof.** In inequality $\sum x^2 - \sum xy \ge \dfrac{3}{4}\max\left\{(x-y)^2; (y-z)^2; (z-x)^2\right\}$, we take $y = [x], z = \{x\}$.

**Theorem 24.** If $x > 0$, then $e^{\{x\}} + e^{[x]} \ge 2 + x$.

**Proof.** In inequality $e^y + e^z \ge 2 + y + z$, we take $y = [x], z = \{x\}$.

**Theorem 25.** If $x \in \mathbb{R}$, then $|\sin[x]| + |\sin\{x\}| + |\cos x| \ge 1$.

**Proof.** In inequality $|\sin a| + |\sin b| + |\cos(a+b)| \ge 1$ we take $a = x, \ b = \{x\}$.

**Theorem 26.** If $x > 0$, then $\left(3[x]^2 + 3[x]\{x\} + \{x\}^2\right) \cdot \left([x]^2 + [x]\{x\} + \{x\}^2\right) \cdot$
$\cdot \left(3\{x\}^2 + 3\{x\}[x] + [x]^2\right) \geq \left(x^2 + [x]\{x\}\right)^3.$

**Proof.** In inequality
$\left(a^2 + ab + b^2\right) \cdot \left(b^2 + bc + c^2\right) \cdot \left(c^2 + ca + a^2\right) \geq (ab + bc + ca)^3$, we take
$a = x, \ b = [x], \ c = \{x\}.$

**Theorem 27.** If $x > 0$, then
$$\frac{[x]}{(3[x] + 2\{x\})([x] + 2\{x\})} + \frac{\{x\}}{(3\{x\} + 2[x])(\{x\} + 2[x])} \geq \frac{11}{48x}.$$

**Proof.** In inequality $\left(\sum x\right) \sum \frac{x}{(2x + y + z)(y + z)} \geq \frac{9}{8}$, we take $y = [x], z = \{x\}$.

**Theorem 28.** If $x > 0$, then $\dfrac{[x]^2}{x + [x]} + \dfrac{\{x\}^2}{x + \{x\}} \geq \dfrac{x\left(2x^2 + 3[x]\{x\}\right)}{(x + [x])(x + \{x\})}.$

**Proof.** In inequality $\sum \dfrac{x^2}{(x + y)(x + z)} \geq \dfrac{3}{4}$, we take $y = [x], z = \{x\}$.

**Theorem 29.** If $x > 1$, then
$$\sqrt{1 + \frac{2\{x\}}{[x]}} + \sqrt{1 + \frac{2[x]}{\{x\}}} \geq 1 + 2\left(\sqrt{\frac{[x]}{x + \{x\}}} + \sqrt{\frac{\{x\}}{x + [x]}}\right).$$

**Proof.** In inequality $\sum \sqrt{\dfrac{y + z}{x}} \leq \sqrt{2} \sum \sqrt{\dfrac{x}{y + z}}$ we take $y = [x], z = \{x\}$.

**Theorem 30.** If $x > 1$, then $\dfrac{\{x\}}{[x]} + \dfrac{[x]}{\{x\}} \geq 1 + 2\left(\dfrac{[x]}{x + \{x\}} + \dfrac{\{x\}}{x + [x]}\right).$

**Proof.** In inequality $\sum \dfrac{y + z}{x} \geq 4 \sum \dfrac{x}{y + z}$, we take $y = [x], z = \{x\}$.

**Theorem 31.** If $x > 0$, then
1). $\min\left\{\left(\sqrt{2} + 1\right)\sqrt{x} + \sqrt{\{x\}}; \left(\sqrt{2} + 1\right)\sqrt{x} + \sqrt{x + [x]}\right\} \geq \sqrt{5([x] + 2\{x\})}$
2). $\left(\sqrt{2} + 1\right)\sqrt{x} + \sqrt{x + \{x\}} \geq \sqrt{5(\{x\} + 2[x])}.$

**Proof.** In $\sqrt{a + b + c} + \sqrt{b + c} + \sqrt{c} \geq \sqrt{a + 4b + 9c}$, we take $a = x, \ b = [x], c = \{x\}$, etc.

**Theorem 32.** If $x \in \mathbb{R}$, then
1). $|\sin x| \leq |\sin[x]| + |\sin\{x\}|$
2). $|\cos x| \leq |\cos[x]| + |\cos\{x\}|$

**Proof.** In inequalities $|\sin(a+b)| \leq |\sin a| + |\sin b|$ and $|\cos(a+b)| \leq |\cos a| + |\cos b|$, we take $a = x$, $b = [x]$.

**Theorem 33.** If $x > 1$, then $6 + \dfrac{\{x\}}{[x]} + \dfrac{[x]}{\{x\}} \geq \left(\sqrt[3]{\dfrac{x}{[x]}} + \sqrt[3]{\dfrac{[x]}{\{x\}}} + \sqrt[3]{\dfrac{\{x\}}{x}}\right)^3$.

**Proof.** In inequality $3\left(\sum a\right)\left(\sum \dfrac{1}{a}\right) \geq \left(\sum \sqrt[3]{\dfrac{a}{b}}\right)$, we take $a = x$, $b = [x]$, $c = \{x\}$.

**Theorem 34.** If $x > 0$, then $\dfrac{[x]}{(x+\{x\})^2} + \dfrac{\{x\}}{(x+[x])^2} \geq \dfrac{1}{8x}$.

**Proof.** In inequality $\sum \dfrac{a}{(b+c)^2} \geq \dfrac{9}{4\sum a}$, we take $a = x$, $b = [x]$, $c = \{x\}$.

**Theorem 35.** If $x > 0$, then
$$\dfrac{x[x]}{(x+\{x\})(2x+[x])} + \dfrac{\{x\}(x+\{x\})}{(x+[x])(2x+\{x\})} \geq \dfrac{[x]+5\{x\}}{12x}.$$
**Proof.** In inequality $\sum \dfrac{a(a+b)}{(b+c)(2a+b+c)} \geq \dfrac{3}{4}$, we take $a = x$, $b = [x]$, $c = \{x\}$.

**Theorem 36.** If $x > 0$, then $\dfrac{[x]}{2x+\{x\}} + \dfrac{\{x\}}{2x+[x]} \leq \dfrac{3[x]^2 + 4[x]\{x\} + 3\{x\}^2}{6x}$.

**Proof.** In inequality $\sum \dfrac{ab}{a+b+2c} \leq \dfrac{1}{4}\sum a$, we take $a = x$, $b = [x]$, $c = \{x\}$.

**Theorem 37.** If $x > 0$, then $\left([x]^5 - [x]^2 + 3\right)\left(\{x\}^5 - \{x\}^2 + 3\right) \geq \dfrac{8x^3}{x^5 - x^2 + 3}$.

**Proof.** In inequality $\prod\left(a^5 - a^2 + 3\right) \geq \left(\sum a\right)^3$, we take $a = x$, $b = [x]$, $c = \{x\}$.

**Theorem 38.** If $x > 0$, then $\dfrac{(2x+[x])^2}{2[x]^2 + (x+\{x\})^2} + \dfrac{(2x+\{x\})^2}{2\{x\}^2 + (x+[x])^2} \leq 5$.

**Proof.** In inequality $\sum \dfrac{(2a+b+c)^2}{2a^2+(b+c)^2} \leq 8$, we take $a = x$, $b = [x]$, $c = \{x\}$.

**Theorem 39.** If $x > 0$, then

1). $\left(x + \sqrt{x[x]} + \sqrt[3]{x[x]\{x\}}\right)^3 \leq 9x^2(x+[x])$

2). $\left([x] + \sqrt{[x]\{x\}} + \sqrt[3]{x[x]\{x\}}\right)^3 \leq 9x^2[x]$

3). $\left(\{x\} + \sqrt{x\{x\}} + \sqrt[3]{x[x]\{x\}}\right)^3 \leq 9x^2(x+\{x\})$.

**Proof.** In inequality $\dfrac{a + \sqrt{ab} + \sqrt[3]{abx}}{2} \leq \sqrt[3]{a\left(\dfrac{a+b}{2}\right)\left(\dfrac{a+b+c}{b}\right)}$, we take $a = x$, $b = [x]$, $c = \{x\}$, etc.

**Theorem 40.** If $x > 0$, then $7(x+[x])^4 + 7(x+\{x\})^4 \geq 3x^4 + 4\left([x]^4 + \{x\}^4\right)$.

**Proof.** In inequality $\sum (a+b)^4 \geq \dfrac{4}{7}\sum a^4$, we take $a = x$, $b = [x]$, $c = \{x\}$.

**Theorem 41.** If $x > 0$, then $\dfrac{\{x\}^2}{(x+\{x\})^2 + [x]^2} + \dfrac{[x]^2}{(x+[x])^2 + \{x\}^2} \geq \dfrac{3}{20}$.

**Proof.** In inequality $\sum \dfrac{(b+c-a)^2}{(b+c)^2 + a^2} \geq \dfrac{3}{5}$, we take $a = x$, $b = [x]$, $c = \{x\}$.

**Theorem 42.** If $x > 0$, then $\sqrt{\dfrac{2x}{x+[x]}} + \sqrt{\dfrac{2[x]}{x}} + \sqrt{\dfrac{2\{x\}}{x+\{x\}}} \leq 3$.

**Proof.** In inequality $\sum \dfrac{2a}{a+b} \leq 3$, we take $a = x$, $b = [x]$, $c = \{x\}$.

**Theorem 43.** If $x > 0$, then $\dfrac{1}{(x+[x])^2} + \dfrac{1}{(x+\{x\})^2 + \{x\}^2} \geq \dfrac{5x^2 - 4[x]\{x\}}{4x^2\left(x^2 + [x]\{x\}\right)}$.

**Proof.** In inequality $\left(\sum xy\right)\left(\sum \dfrac{1}{(x+y)^2}\right) \geq \dfrac{9}{4}$, we take $y = [x]$, $z = \{x\}$.

**REFERENCES:**

[1] Mihály Bencze, Inequalities (manuscript), 1982.
[2] Collection of "Octogon" Mathematical Magazine (1993-2006).